\documentclass{elsart3-1}

\usepackage{amsmath,  amsopn, amsfonts}

\usepackage{amssymb}

\usepackage[english,francais]{babel}

\newtheorem{e-proposition}[theorem]{Proposition}

\newtheorem{e-definition}[theorem]{Definition\rm}
\newtheorem{remark}{\it Remark\/}

\renewcommand{\theequation}{\arabic{equation}}
\def\og{\leavevmode\raise.3ex\hbox{$\scriptscriptstyle\langle\!\langle$~}}
\def\fg{\leavevmode\raise.3ex\hbox{~$\!\scriptscriptstyle\,\rangle\!\rangle$}}

\def\11{{\rm 1~\hspace{-1.4ex}l} }

\begin{document}

\begin{frontmatter}



\selectlanguage{english}
\title{Dispersion for the wave and the Schr\"{o}dinger equations outside strictly convex obstacles and counterexamples}

\vspace{-2.6cm}

\selectlanguage{francais}
\title{Estimations de dispersion pour l'\'equation des ondes et de Schr\"odinger \`a l'ext\'erieur des obstacles strictement convexes }


 \author{Oana Ivanovici}
    \ead{oana.ivanovici@unice.fr}
\ead[url]{http://www.math.unice.fr/\~{}ioana}
  \address{ CNRS et Universit\'e C\^ote d'Azur\\
  Laboratoire J. A. Dieudonn\'e, UMR CNRS 7351\\
  Parc Valrose, 06108 Nice Cedex 02,
    France}

  \author{Gilles Lebeau}
  \ead{gilles.lebeau@unice.fr}
  
  \address{Universit\'e C\^ote d'Azur\\
  Laboratoire J. A. Dieudonn\'e, UMR CNRS 7351\\
    Parc Valrose,
    06108 Nice Cedex 02,
    France}
    
        \thanks{The authors were partially supported by ERC project SCAPDE. The authors would like to thank Centro di Giorgi, Pisa for the warm welcome during the summer 2015 when this article has started.
    } 
\selectlanguage{english}

\begin{abstract}
The purpose of this note is to prove dispersive estimates for the wave and the Schr\"odinger equations outside strictly convex obstacles in $\mathbb{R}^d$. If $d=3$, we show that for both equations, the linear flow satisfies the (corresponding) dispersive estimates as in $\mathbb{R}^3$. In higher dimensions $d\geq 4$ and if the domain is the exterior of a ball in $\mathbb{R}^d$, we show that losses in dispersion do appear and this happens at the Poisson spot.

\vskip 0.5\baselineskip
\selectlanguage{francais}
\noindent{\bf R\'esum\'e}
\vskip 0.5\baselineskip
\noindent
L'objet de cette note est de d\'emontrer des estimations de dispersion pour l'\'equation des ondes et de Schr\"odinger \`a l'ext\'erieur d'un obstacle strictement convexe de $\mathbb{R}^d$. Si $d=3$, on d\'emontre  que, pour chacune des deux \'equations, le flot lin\'eaire v\'erifie les estimations de dispersion comme dans $\mathbb{R}^3$. En dimension plus grande $d\geq 4$, on d\'emontre que des pertes dans la dispersion apparaissent \`a l'ext\'erieur d'une boule de $\mathbb{R}^d$ et cela arrive au point de Poisson.

\end{abstract}
\end{frontmatter}

\selectlanguage{francais}
\section*{Version fran\c{c}aise abr\'eg\'ee}

Pour l'\'equation des ondes, dans le cas Euclidien, la forme explicite du flot permet d'obtenir les estimations de dispersion 
\[
\|\chi(hD_t)e^{\pm it\sqrt{-\Delta_{\mathbb{R}^d}}}\|_{L^1(\mathbb{R}^d)\rightarrow L^{\infty}(\mathbb{R}^d)}\leq C(d)h^{-d}\min\{1,(h/|t|)^{\frac{d-1}{2}}\}, \quad \chi\in C_{0}^\infty
(]0,\infty[).
\]
Pour l'\'equation de Schr\"odinger, les estimations de dispersion s'\'ennoncent comme suit:
\[
\|e^{\pm it\Delta_{\mathbb{R}^d}}\|_{L^1(\mathbb{R}^d)\rightarrow L^{\infty}(\mathbb{R}^d)}\leq C(d)|t|^{-d/2}.
\]

Notre but est d'obtenir des estimations de dispersion \`a l'ext\'erieur d'un obstacle strictement convexe. Plusieurs r\'esultats positifs sur les effets dispersifs ont \'et\'e obtenus r\'ecemment dans ce contexte: cependant, la question de savoir si les estimations de dispersion \'etaient vraies ou non est rest\'ee ouverte, m\^eme \`a l'ext\'erieur d'une boule. Puisqu'il n'y a pas de concentration apparente d'\'energie, comme ledans le cas d'un domaine non-captant quelconque (pour lequel les portions concaves du bord peuvent agir comme des miroirs et re-focaliser les paquets d'ondes) on pourrait raisonnablement penser que les estimations de dispersion devraient \^etre v\'erifi\'ees \`a l'ext\'erieur d'un convexe (voir l'ext\'erieur d'une ball \cite{lismza12} dans le cas des fonctions \`a sym\'etrie spherique). On montre ici que c'est effectivement le cas en dimension $3$, par contre en dimension plus grande on construit des contre-exemples explicites \`a l'ext\'erieur d'une boule.

\begin{thm}
Soit $\Theta\subset\mathbb{R}^3$ un domaine compact avec bord r\'egulier, strictement convexe et soit $\Omega=\mathbb{R}^3\setminus \Theta$. Soit $\Delta$ le Laplacien dans $\Omega$ avec condition de Dirichlet au bord. Alors 
\begin{enumerate}
\item les estimations de dispersion pour le propagateur des ondes dans $\Omega$ sont v\'erifi\'ees comme dans $\mathbb{R}^3$:
\[
\|\chi(hD_t)e^{\pm it\sqrt{-\Delta}}\|_{L^1(\Omega)\rightarrow L^{\infty}(\Omega)}\leq Ch^{-3}\min\{1,\frac{h}{|t|}\}.
\]
\item  les estimations de dispersion pour le flot de Schr\"odinger dans $\Omega$ sont v\'erifi\'ees comme dans $\mathbb{R}^3$.

\end{enumerate}
\end{thm}
On remarque qu'une perte dans la dispersion pourrait \^etre li\'ee (de fa\c con informelle) \`a la pr\'esence d'un point de concentration : ces points apparaissent lorsque des rayons optiques (envoy\'es d'une m\^eme source dans des directions diff\'erentes) cessent de diverger. Le principe de Huygens \'enonce que lorsque la lumi\`ere \'eclaire un obstacle circulaire, chaque point de l'obstacle se comporte \`a son tour comme une nouvelle source lumineuse ponctuelle; tous les rayons lumineux issus des points de la circonf\'erence de l'obstacle se concentrent au centre de l'ombre et d\'ecrivent le m\^eme chemin optique; il en r\'esulte une tache lumineuse au centre de l'ombre (le point de Poisson).
Par cons\'equent, l'intuition nous dit que s'il y a une perte dans la dispersion, elle devrait appara\^itre au point de Poisson.

\begin{thm}
Pour $d\geq 4$ on pose $\Theta=B_d(0,1)$ la boule unit\'e de $\mathbb{R}^d$. Soit $\Omega_d=\mathbb{R}^d\setminus B_d(0,1)$ et soit $\Delta_d$ le Laplacien dans $\Omega_d$ avec condition de Dirichlet. Au point de Poisson, les estimations de dispersion pr\'ec\'edentes (o\`u $\mathbb{R}^d$ et $\Delta_{\mathbb{R}^d}$ sont remplac\'es par $\Omega_d$ et $\Delta_d$) ne sont plus v\'erifi\'ees. Precis\'ement, soient $Q_{\pm}(r)$ les points source et d'observation situ\'es \`a distance $r$ du centre $O$ de la boule $B_d(0,1)$, sym\'etriques par rapport \`a $O$;  alors, si $r=\gamma h^{-1/3}$, avec $\gamma$ \`a valeurs dans un compact de $(0,\infty)$,
\begin{itemize}
\item pour le propagateur des ondes et pour $t\simeq 2\gamma h^{-1/3}$
\[
\Big|(\chi(hD_t)e^{it\sqrt{|\Delta|}}(\delta_{Q_+(\gamma h^{-1/3})})\Big| (Q_-(\gamma h^{-1/3}))\simeq h^{-d}\Big(\frac ht\Big)^{\frac{d-1}{2}} h^{-\frac{d-3}{3}}, 
\]
\item pour le propagateur de Schr\"odinger classique et pour $t\simeq h^{1/3}$
\[
\Big|(\chi(hD_t)e^{it\Delta}(\delta_{Q_+(\gamma h^{-1/6})})\Big| (Q_-(\gamma h^{-1/6}))\simeq h^{-1-\frac d6-\frac{d-3}{6}}.
\]
\end{itemize}
Pour $d\geq 4$, ces estimations contredissent les estimations du cas plat $\mathbb{R}^d$.
\end{thm}
\selectlanguage{english}
\section{Introduction}
\label{}



\selectlanguage{english}
 
In the Euclidean case, the explicit form of the wave propagator yields the following dispersive estimate
\begin{equation}\label{disprd}
\|\chi(hD_t)e^{\pm it\sqrt{-\Delta_{\mathbb{R}^d}}}\|_{L^1(\mathbb{R}^d)\rightarrow L^{\infty}(\mathbb{R}^d)}\leq C(d)h^{-d}\min\{1,(h/|t|)^{\frac{d-1}{2}}\}, \quad \chi\in C_{0}^\infty
(]0,\infty[).
\end{equation}
Concerning the Schr\"odinger equation, the dispersive estimates read as follows:

\begin{equation}\label{dispschrodrd}
\|e^{\pm it\Delta_{\mathbb{R}^d}}\|_{L^1(\mathbb{R}^d)\rightarrow L^{\infty}(\mathbb{R}^d)}\leq C(d)|t|^{-d/2}.
\end{equation}

Our aim in the present paper is to obtain dispersive estimates outside strictly convex obstacles. While many positive results on dispersive effects had been established lately in this context, the question about whether or not dispersion did hold remained open, even for the exterior of a ball. Since there is no apparent concentration of energy, like in the case of a generic non-trapping obstacle (where concave portions of the boundary can act as mirrors and refocus wave packets), one would expect dispersive estimates to hold outside strictly convex obstacles (see the exterior of a ball \cite{lismza12} for spherically symmetric functions). We prove that this is indeed the case in dimension three, while in higher dimensions we provide explicit counterexamples for the exterior of a ball. 

\begin{thm}\label{thmdisp3D}
Let $\Theta\subset\mathbb{R}^3$ be a compact domain with smooth, strictly convexe boundary and let $\Omega=\mathbb{R}^3\setminus \Theta$. Let $\Delta$ denote the Dirichlet Laplace operator in $\Omega$. Then 
\begin{enumerate}
\item the dispersive estimates for the wave flow in $\Omega$ do hold like in $\mathbb{R}^3$:
 \begin{equation}\label{dispomega3}
\|\chi(hD_t)e^{\pm it\sqrt{-\Delta}}\|_{L^1(\Omega)\rightarrow L^{\infty}(\Omega)}\leq Ch^{-3}\min\{1,\frac{h}{|t|}\}.
\end{equation}
\item  the dispersive estimates for the classical Schr\"odinger flow in $\Omega$ hold like in $\mathbb{R}^3$.
\end{enumerate}
\end{thm}
We remark that a loss in dispersion may be informally related to a cluster point : such clusters occur because optical rays (sent from the same source along different directions) are no longer diverging from each other.
When light shines on a circular obstacle, Huygens's principle says that every point of the obstacle acts as a new point source of light, so all the light passing close to a perfectly circular object concentrate at the perfect center of the shadow behind it; this results in a bright spot at the shadow's center (the Poisson spot). Therefore, our intuition tells us that if there is a location where dispersion could fail, this will happen at the Poisson spot.

\begin{thm}\label{thmCE}
Let $d\geq 4$ and let $\Theta=B_d(0,1)$ be the unit ball in $\mathbb{R}^d$. Set $\Omega_d=\mathbb{R}^d\setminus B_d(0,1)$ and let $\Delta_d$ denote the Laplace operator in $\Omega_d$. Then at the Poisson spot the dispersive estimates \eqref{disprd},\eqref{dispschrodrd} (with $\mathbb{R}^d$ and $\Delta_{\mathbb{R}^d}$ replaced by $\Omega_d$ and $\Delta_d$) fail. Precisely, let $Q_{\pm}(r)$ be the source and the observation points at (same) distance $r$ from the ball $B_d(0,1)$, symmetric with respect to the center of the ball, then, taking $r=\gamma h^{-1/3}$, with $\gamma$ in a compact subset of $(0,\infty)$ yields
\begin{itemize}
\item for the wave flow and for $t\simeq 2\gamma h^{-1/3}$
\begin{equation}
\Big|(\chi(hD_t)e^{it\sqrt{|\Delta|}}(\delta_{Q_+(\gamma h^{-1/3})})\Big| (Q_-(\gamma h^{-1/3}))\simeq h^{-d}\Big(\frac ht\Big)^{\frac{d-1}{2}} h^{-\frac{d-3}{3}}, 
\end{equation}
\item for the classical Schr\"odinger flow and for $t\simeq h^{1/3}$
\[
\Big|(\chi(hD_t)e^{it\Delta}(\delta_{Q_+(\gamma h^{-1/6})})\Big| (Q_-(\gamma h^{-1/6}))\simeq h^{-1-\frac d6-\frac{d-3}{6}}.
\]
\end{itemize}
For $d\geq 4$, these estimates contradict the usual ones \eqref{disprd}, \eqref{dispschrodrd} in $\mathbb{R}^d$.
\end{thm}

\section{General setting for the wave flow outside a ball in $\mathbb{R}^d$}
In this note we give a sketch of the proof of Theorem \ref{thmdisp3D} only in the case of the wave equation outside a ball. The general case will be dealt with in \cite{IL17}.
Let $\Omega=\mathbb{R}^d\setminus B_d(0,1)$, $\partial\Omega=\mathbb{S}^{d-1}$. Let $Q_0\in \Omega$ and $\delta_{Q_0}$ the Dirac distribution at $Q_0$. Let also $\Delta_{\mathbb{R}^d}$ denote the Laplace operator in the whole space $\mathbb{R}^d$ and $U(t,Q,Q_0)=\cos(t\sqrt{|\Delta_{\mathbb{R}^d}|})(\delta_{Q_0})(Q)$ be the solution to
\[
   (\partial^2_t-
 \Delta_{\mathbb{R}^d}) U=0  \;\; \text{ in } \Omega; 
 U|_{t=0}=\delta_{Q_0}, \; \partial_t U|_{t=0}= 0;
 U|_{\partial\Omega}=0. 
\]
Then $\partial_n U  |_{\partial\Omega}=\partial_n U_{free}|_{\partial\Omega}-{\bf N}(U_{free}|_{\partial\Omega})$,
where $U_{free}$ is the free wave in $\mathbb{R}^d$, $\vec{n}$ is the outward unit normal to $\partial\Omega$ pointing towards $\Omega$ and ${\bf N}$ is the Neumann operator. Define
\[
\underline{U}(t,Q,Q_0):=
U(t,Q,Q_0), \text{ if } Q\in\Omega; \quad 
\underline{U}(t,Q,Q_0):=0, \text{ if } Q\in \overline{B_d(0,1)}.
\]
Then $\underline{U}$ satisfies the following equation:
\[
   (\partial^2_t-
 \Delta_{\mathbb{R}^d}) \underline{U}=\partial_n U|_{\partial\Omega}\otimes \delta_{\partial\Omega} \;\; \text{ in } \mathbb{R}^d;
 \underline{U}|_{t=0}=\delta_{Q_0}, \; \partial_t \underline{U}|_{t=0}=0,
\]
which yields
$ \underline{U}|_{t>0} =\square^{-1}_+\Big(\partial_n U|_{\partial\Omega}\otimes \delta_{\partial\Omega}|_{t>0}\Big)+{ U_{free}|_{t>0}}$, where
\[
\square^{-1}_{+}F(t)=\int_{-\infty}^t R(t-t')*F(t')dt',  \quad \widehat{R(t,\xi)}=\frac{\sin (t|\xi|)}{|\xi|}.
\]
\section{ Sketch of proof of Theorem \ref{thmdisp3D}}
In dimension $d=3$, from the last formula and the form of $\partial_n U|_{\partial\Omega}$ in terms of $U_{free}$ and ${\bf N}$ we find
\begin{multline}\label{toestimate}
\square^{-1}_+\Big(\Big(\partial_n U_{free}|_{\partial\Omega}-N(U_{free}|_{\partial\Omega})\Big)\otimes \delta_{\partial\Omega}|_{t>0}\Big)
=\frac{1}{4\pi}\int_{\partial\Omega}\frac{(\partial_n U_{free}|_{\partial\Omega}-{\bf N}(U_{free}|_{\partial\Omega}))(t-|Q-P|,P)}{4\pi |P-Q|}d\sigma(P).
\end{multline}
In order to prove Theorem \ref{thmdisp3D} (in dimension $3$) we are reduced to obtaining bounds for \eqref{toestimate}. For that, we use the Melrose and Taylor parametrix which provides the form of the solution near the glancing regime in terms of Airy function. Outside a neighborhood of the glancing region it is easy to see that the dispersive estimates hold true. The next theorem is due to Melrose and Taylor and holds for $\partial\Omega$ strictly concave :
\begin{thm} $\exists \theta,\zeta$ phase functions near the glancing region, $\exists$ $a, b$ symbols (with $a$ elliptic, $b|_{\partial\Omega}=0$) such that, if $V$ is a solution in $\Omega$ to
\[
 (\tau^2+\Delta )V\in O_{C^{\infty}}(\tau^{-\infty}),
\]
then there exists $F$ such that
\[
V(\tau,Q)=\frac{\tau}{2\pi}\int e^{i\tau\theta(Q,\eta)} \Big(aA+\tau^{-1/3}b A'\Big)(\tau^{2/3}\zeta(Q,\eta))\hat{F}(\tau\eta)d\eta,
\]
where we set $A(z)=Ai(-z)$, where $Ai$ is the Airy function which satisfies $Ai''(z)=zAi(z)$.
\end{thm}
For $d=3$, we take ${V=\widehat{U_{free}}(\tau,P,Q_0)=\frac{\tau}{|P-Q_0|}e^{-i\tau |P-Q_0|}}$.
We introduce polar coordinates: since $\partial\Omega=\mathbb{S}^2$, a point in $\Omega$ can be written as $(r,\varphi,\omega)$, $r>1$, $\varphi\in [0,\pi]$, $\omega\in \mathbb{S}^1$. We can always assume that the source point $Q_0$ has coordinates $(r_0,0,.)$, $r_0>1$. We define the apparent contour $\mathcal{C}_Q$ of a point $Q\in \Omega$ as the boundary of the set of points that can be "viewed" from $Q$. Therefore, $\mathcal{C}_{Q_0}=\{(1,\varphi,\omega)|\quad \cos(\varphi)=1/r_0\}$. 
\begin{remark}
When $\Omega=\mathbb{R}^d\setminus B_d(0,1)$, $d\geq 2$, the functions $\theta$ and $\zeta$ can be taken under the following form 
\[
\theta(\varphi,\eta)=\varphi\eta,\quad
\zeta(r,\eta)=\eta^{2/3}l(\frac{r}{\eta}),
\]
where $\eta\simeq 1$ and $l(z)=(z-1)(1+ O(z-1)) \text{ for } z \text{ close to } 1$. Then $\zeta_0(\eta):=\zeta(1,\eta)=\eta^{2/3}l(1/\eta)$.
\end{remark}
Let $P=(1,\varphi,\omega)\in \mathbb{S}^2$ and $a_0=a|_{\partial\Omega}$, then the trace of the free wave on the boundary $\mathbb{S}^2$ reads as
\begin{equation}\label{formhatF}
\widehat{U_{free}}(\tau,P,Q_0)=\frac{\tau}{|P-Q_0|}e^{-i\tau |P-Q_0|}
=\frac{\tau}{2\pi}\int e^{i\tau \varphi\eta} a_0 A(\tau^{2/3}\zeta_0(\eta))\hat{F}(\tau\eta)d\eta.
\end{equation}
Let $A_+(z)=A(e^{i\pi/3}z)$ : then $A_+(z)$ doesn't vanish for $z\in \mathbb{R}$. We compute
\[
\partial_n\widehat{U_{free}}(\tau,P,Q_0)=\frac{\tau}{2\pi}\int e^{i\tau \varphi\eta}\Big(\tilde aA+\tau^{2/3}\tilde b A'\Big)(\tau^{2/3}\zeta_0(\eta))\hat{F}(\tau\eta)d\eta,
\]
\[
\widehat{{\bf N}(U_{free}|_{\mathbb{S}^2})}(\tau,P,Q_0)=\frac{\tau}{2\pi}\int e^{i\tau \varphi\eta}\Big(\tilde{a}A+\tau^{2/3}\tilde{b}\frac{A'_+}{A_+}A\Big)(\tau^{2/3}\zeta_0(\eta))\hat{F}(\tau\eta)d\eta.
\]
We obtain an explicit form for $\square^{-1}_+\Big(\partial_n U|_{\partial\Omega}\otimes \delta_{\partial\Omega}|_{t>0}\Big)$ as follows
\[
\square^{-1}_+\Big(\Big(\partial_n U_{free}|_{\mathbb{S}^2}-{\bf N}(U_{free}|_{\mathbb{S}^2})\Big)\otimes \delta_{\mathbb{S}^2}\Big)=\int e^{i\tau t}\chi(h \tau)({I_F+I_D+I_R})(\tau,Q,Q_0)d\tau,
\]
where 
\begin{itemize}
\item ${ I_F}$ is the direct wave: the phase is the phase of the free wave and the amplitude is just the amplitude of the free wave cutoff near the shadow ($\Rightarrow $ {OK} for dispersion);  
\item ${I_D}$ is the diffracted wave: it corresponds to a neighborhood of $\mathcal{C}_{Q_0}$ on the boundary of size $\tau^{-1/3}$ in $\varphi$ and of size $\tau^{-2/3}$ in angle around the glancing direction (around $\eta=1$).
\item ${I_R}$ is the reflected wave: the phase has a singular, Airy type term (easy to deal with).
\end{itemize}
Since difficulties appear near rays issued from $Q_0$ which hit the boundary without being deviated, only the diffracted wave  part (containing $I_D$) will be dealt with here. Notice that this is the regime which provides counter-examples in higher dimensions. We have
\[
I_D(\tau,Q,Q_0)=\int_{P=(1,\varphi,\omega)\in\mathbb{S}^2}\frac{1}{|P-Q|}e^{-i\tau|P-Q|}\frac{\tau}{2\pi}\int e^{i\tau \varphi \eta}\tau^{2/3}\tilde{b}
\frac{\chi_0(\tau^{2/3}\zeta_0(\eta))}{(AA_+)(\tau^{2/3}\zeta_0(\eta))}\widehat{F}(\tau\eta)d\eta d\varphi d\omega,
\]
where $\chi_0$ is supported near $0$, $\eta\simeq 1$ and $\tilde b$ is an elliptic symbol. 

The  phase function of $I_D$ equals ($-|P-Q|+\varphi\eta+$ phase of $\hat{F}$) 
and reads as $-|P-Q|+\varphi\eta-\varphi_0\eta-|P_0-Q_0|$,
where $P_0\in \mathcal{C}_{Q_0}$ has coordinates $P_0=(1,\varphi_0,.)$, $\cos(\varphi_0)=1/r_0$.
The symbol of $I_D$ is of the form $\frac{\tau}{|P-Q|}\tau^{2/3}\times\chi_0\times (\tau^{-1/3} \frac{\tau}{|P_0-Q_0|})$, where the factor in brackets comes from $\hat{F}$, obtained from \eqref{formhatF} as an oscillatory integral with critical points of order precisely $2$ on $\mathcal{C}_{Q_0}$. It will be enough to prove that $|I_D(\tau,Q,Q_0)|\lesssim \frac{\tau}{t}$. 

Let $Q=(r_Q,\varphi_Q,\omega_Q)$ be an observation point in $\Omega$; in $I_D$, the only dependence in $\omega$ comes from $|P-Q|$ since
\[
 |P-Q|=\Big(1+r^2_Q-2r_Q\sin(\varphi)\sin(\varphi_Q)\cos(\omega-\omega_Q)-2r_Q\cos(\varphi)\cos(\varphi_Q)\Big)^{1/2}.
\]
The critical points with respect to $\omega$ satisfy
$\partial_{\omega} |P-Q|=\frac{r_Q}{|P-Q|}\sin(\varphi)\sin(\varphi_Q)\sin (\omega-\omega_Q)$.
\begin{itemize}
\item If $\sin(\varphi_Q)\neq 0$ $\Rightarrow $ the critical points satisfy $\sin(\omega-\omega_Q)=0$ which yields $\omega=\omega_Q$ or $\omega=\omega_Q+\pi$; this means that the stationary points $P$ on the boundary belong to a circle situated in the plane $Q_0-O-Q$; all the computations are explicit and provide the announced result;

\item If $\sin(\varphi_Q)=0$ the derivative vanishes everywhere. In this case the points $Q$, $O$ and $Q_0$ are colinear and the integration in $\omega$ does not provide negative factors of $\tau$. It is easy to see that the integration with respect to $\varphi$ provides a power of $\tau^{-1/3}$ (corresponding to the critical points $\varphi$ such that $(0,\varphi,\omega)\in\mathcal{C}_Q$ which are degenerate of order $2$) and the integration with respect to $\eta$ provides a factor $\tau^{-2/3}$, due to the localisation $\chi_0$. It remains to show that the remaining factor  $\frac{\tau^{1/3}}{|P-Q| |P_0-Q_0|}$ must be bounded by $\frac Ct$, for some $C>0$; indeed, $\mathcal{C}_{Q_0}$ is defined by $\cos(\varphi_0)=\frac{1}{r_0}$, while $\mathcal{C}_Q$ is defined by  $\cos(\varphi)=\frac{1}{r_Q}$, and $|\varphi-\varphi_0|\lesssim \tau^{-1/3}$ (since otherwise, by integrations by parts with respect to $\eta$, we get a $O(\tau^{-\infty})$ contribution). Notice that we have also used that $t\simeq |Q_0-P_0|+|P-Q|$.
\end{itemize}
This allows to achieve the proof of Theorem \ref{thmdisp3D}.

\section{Sketch of proof of Theorem \ref{thmCE} for the wave flow outside a ball in $\mathbb{R}^d$, $d\geq 4$}
Let $S$, $N$ denote the south pole and the north pole, respectively. Let $Q_+(r)$ be a point on $OS$ axis, at distance $r$ from $O$ and let $Q_-(r)$ denote its symmetric with respect to $O$ on the $ON$ axis, where $O\in\mathbb{R}^{d}$ is the centre of the ball. We let $r_0=\gamma h^{-1/3}$ for some $\gamma$ in a compact set of $(0,\infty)$ and let $Q_0=Q_+(r_0)$ and $Q=Q_-(r_0)$. The counterexample to dispersion  comes from the diffracted part, which, in this case, takes the form 
\[
\int e^{it\tau}\chi(h\tau) I_D(\tau,Q_-,Q_+)d\tau,
\]
where for $d\geq 3$
\begin{multline*}
I_D(\tau,Q,Q_0)=\int_{P=(1,\varphi,.)\in\mathbb{S}^{d-1}}\frac{\tau^{\frac{d-3}{2}}}{|P-Q|^{\frac{d-1}{2}}}e^{-i\tau|P-Q|}\Sigma_d(|P-Q|\tau)\\\times \frac{\tau}{2\pi}\int_{\mathbb{R}} e^{i\tau \varphi \eta}\tau^{2/3}\tilde{b}
\frac{\chi_0(\tau^{2/3}\zeta_0(\eta))}{(AA_+)(\tau^{2/3}\zeta_0(\eta))}\widehat{F}(\tau\eta)d\eta d\varphi,
\end{multline*}
where $\Sigma_d$ is a symbol of degree $0$ which satisfies
\[
\widehat{U}_{free}(\tau,P,Q_0)=\frac{\tau^{\frac{d-1}{2}}}{|P-Q_0|^{\frac{d-1}{2}}}e^{-i\tau |P-Q_0|}\Sigma_d(|P-Q_0|\tau),\quad \text{ for } |P-Q_0|\gg \tau^{-1}.
\]
With $P_0=(1,\varphi_0,.)$, $\cos(\varphi_0)=1/r_0$, on the apparent contour $\mathcal{C}_{Q_0}$, we have
\[
\widehat{F}(\tau\eta)=\tau^{-1/3}\frac{\tau^{\frac{d-1}{2}}}{|P_0-Q_0|^{\frac{d-1}{2}}}e^{i\tau\eta\varphi_0}e^{-i\tau |P_0-Q_0|}\Sigma (|P_0-Q_0|,\tau),
\]
where $\Sigma$ is obtained from $\Sigma_d$ after applying the stationary phase with degenerate critical point $P_0$ on $\mathcal{C}_{Q_0}$. Notice that the observation point $Q=Q_-(r_0)$ is such that $\sin(\varphi_Q)=0$, since $Q_+$, $O$ and $Q_-$ are on the same line and, due to rotational symmetry, in the integral defining $I_D$ the phase function does not depend on $\omega$. We obtain by explicit computations 
\[
\Big|\int e^{it\tau}\chi(h\tau) I_D(\tau,Q,Q_0)d\tau\Big|\simeq C \frac1h \frac{h^{-(d-2+\frac 13)}}{|P_0-Q_0|^{d-1}}.
\]
Since we must have $t\simeq 2\gamma h^{-1/3}$, it follows that
\[
\Big|\int e^{it\tau}\chi(h\tau) I_D(\tau,Q,Q_0)d\tau\Big|\simeq \frac{C}{\gamma^{\frac{d-1}{2}}} h^{-d} \frac{h^{\frac{d-1}{2}}}{|t|^{\frac{d-1}{2}}} h^{{-\frac{d-3}{3}}}.
\]
For $d=3$ this coincide with the usual estimates \eqref{disprd} of $\mathbb{R}^3$. However, for $d\geq 4$ there is a loss coming from the factor $h^{-\frac{d-3}{3}}$ for $\gamma$ in a fixed compact of $(0,\infty)$.

\end{document}